\documentclass[%
 aip,cha,
amsmath,amssymb,
reprint,%
]{revtex4-1}
\usepackage{graphicx}% Include figure files
\usepackage{hyperref}
\usepackage{amsmath}  
\usepackage{amssymb}
\usepackage{wrapfig}
\usepackage{gensymb}
\usepackage{enumitem}
\usepackage{indentfirst}
\usepackage{wrapfig}
\usepackage{graphicx}
\usepackage{subcaption}
\usepackage{tikz}
\usetikzlibrary{arrows}
\usetikzlibrary{matrix}
\usepackage{xcolor}
\usepackage{bm}
\usepackage{subcaption}
\usepackage[capitalise]{cleveref}

\newcommand{\ts}{\textsuperscript}

\newtheorem{example}{Example}[section]
\newtheorem{theorem}{Theorem}[section]
\newtheorem{corollary}{Corollary}[theorem]

\tikzset{
  pics/carc/.style args={#1:#2:#3}{
    code={
      \draw[pic actions] (#1:#3) arc(#1:#2:#3);
    }
  }
}

\begin{document}

\author{Anastasiya Salova}
\email[]{avsalova@ucdavis.edu}
%\thanks{}
\affiliation{Department of Physics, University of California, Davis, CA 95616, USA}
\affiliation{Complexity Sciences Center, University of California, Davis, CA 95616, USA}

\author{Jeffrey Emenheiser}
%\thanks{}
\affiliation{Department of Physics, University of California, Davis, CA 95616, USA}
\affiliation{Complexity Sciences Center, University of California, Davis, CA 95616, USA}

\author{Adam Rupe}
\affiliation{Department of Physics, University of California, Davis, CA 95616, USA}
\affiliation{Complexity Sciences Center, University of California, Davis, CA 95616, USA}
\author{James P. Crutchfield}
\affiliation{Department of Physics, University of California, Davis, CA 95616, USA}
\affiliation{Complexity Sciences Center, University of California, Davis, CA 95616, USA}
\affiliation{Santa Fe Institute, Santa Fe, NM 87501, USA}

\author{Raissa M. D'Souza}
%\thanks{}
\affiliation{Complexity Sciences Center, University of California, Davis, CA 95616, USA}
\affiliation{Santa Fe Institute, Santa Fe, NM 87501, USA}
\affiliation{
Department of Computer Science, University of California, Davis, CA 95616, USA}
\affiliation{Department of Mechanical and Aerospace Engineering, University of California, Davis, CA 95616, USA}
%\preprint{}

\title{Koopman Operator and its Approximations for Systems with Symmetries}




\date{\today}

\begin{abstract}
Nonlinear dynamical systems with symmetries exhibit a rich variety of behaviors, including complex attractor-basin portraits and enhanced and suppressed bifurcations. Symmetry arguments provide a way to study these collective behaviors and to simplify their analysis. The Koopman operator is an infinite dimensional linear operator that fully captures a system's nonlinear dynamics through the linear evolution of functions of the state space. Importantly, in contrast with local linearization, it preserves a system's global nonlinear features. We demonstrate how the presence of symmetries affects the Koopman operator structure and its spectral properties. \textcolor{black}{In fact, we show that symmetry considerations can also simplify finding the Koopman operator approximations using the extended and kernel dynamic mode decomposition methods (EDMD and kernel DMD).} Specifically, representation theory allows us to demonstrate that an isotypic component basis induces block diagonal structure in operator approximations, revealing hidden organization. Practically, if the data is symmetric, the EDMD and kernel DMD methods can be modified to give more efficient computation of the Koopman operator approximation and its eigenvalues, eigenfunctions, and eigenmodes. Rounding out the development, we discuss the effect of measurement noise.
\end{abstract}

%\pacs{}% insert suggested PACS numbers in braces on next line

\maketitle %\maketitle must follow title, authors, abstract and \pacs


\begin{quotation}
Many natural and engineered dynamical systems --- power grid networks and biological regulatory networks, to mention two --- exhibit symmetries in their connectivity structure and in their internal dynamics. Some have time-reversal symmetry, others rotational and spatial translation invariance, and others still, combinations. These symmetries are often key for understanding the behavior of systems. For instance, a well-known example that exploits the connection between system state and interconnection symmetries arises in locomotion, where spatio-temporal symmetries of the observed locomotion behavior constrain the structure of neural circuits that generate these patterns. For network systems in particular, symmetries in the connectivity structure are of fundamental importance.  For instance, the structural symmetries of a network of identical oscillators can determine its admissible patterns of symmetry-breaking. That said, additional information beyond the knowledge of the network structure is often required to address more detailed questions about a system's dynamics, such as whether a particular configuration is stable in a given parameter regime. In these cases, the system's linearization near the steady state can be combined with interconnection symmetry to provide the answer. However, these linearization methods are only valid on a local subset of the state space and therefore are not sufficient for global characteristics of nonlinear dynamical systems, such as their attractors, basins, and transients. The Koopman operator, in contrast, is a linear infinite-dimensional operator that evolves the functions on the state space that is valid on the entire state space. We show how to combine symmetry considerations with the Koopman analysis to study nonlinear dynamical systems with symmetries. We use representation theory to determine the effect of symmetries on the Koopman operator and its approximations, drawing out how local dynamical symmetries interact with symmetries arising from the connectivity of system variables. This, in turn, allows us to modify data-driven Koopman operator approximation algorithms to make them more efficient when applied to dynamical systems with symmetries. We illustrate our findings in a simple network of coupled Duffing oscillators with symmetries in individual oscillator dynamics and in their physical couplings.
\end{quotation}
\section{Introduction}
%\label{}
Symmetries of dynamical systems manifest themselves in asymptotic dynamics, bifurcations, and attractor basin structure. Symmetries play a crucial role in guiding the emergence of synchronization and pattern formation, behaviors broadly observed in natural and engineered systems. Methods from group theory, representation theory, and equivariant bifurcation theory provide useful tools to study the common features of systems with symmetries \cite{golubitsky2003symmetry,
golubitsky2012singularities, golubitsky2015recent}.  

Dynamical elements organized into a network are an important class of dynamical systems that often exhibit these behaviors, especially when symmetries appear in both network structure and the dynamics of the individual nodes. Studying the effect of symmetries in network topology of synthetic and real-life systems using computational group theory is an active area of research \cite{macarthur2008symmetry, pecora2014cluster, pecora2017discovering}. Those symmetries lead to phenomena such as full synchronization, cluster synchronization, and formation of exotic steady states such as chimeras \cite{sorrentino2016complete, cho2017stable, matheny2019exotic, emenheiser2016patterns}. Moreover, topological symmetries underlying cluster synchronization of coupled identical elements assist in analyzing the stability of these fully synchronous cluster states \cite{pecora2014cluster, hart2019topological}. For networks of identical coupled oscillators, the form of their limit cycle solutions and the form of their bifurcations can be derived from symmetry considerations \cite{ashwin1992dynamics}. Symmetries are also key in determining network controllability and observability. For example, Refs. \cite{rubin1972controllability, dellnitz2017sensing} explored the effect of explicit network symmetries for linear time-independent and time-dependent networks. Similarly, Refs. \cite{letellier2002investigating, whalen2015observability} considered nonlinear network motifs with symmetries and studied how the presence of different types of structural symmetries affect the observability and controllability of the system. Ref. \cite{mesbahi2019nonlinear}, similar to our approach, uses the Koopman operator formalism (discussed below). They provide analytic results that link the presence of permutational symmetries in dynamical systems to their observability properties.

Many dynamical systems of current interest are high dimensional and nonlinear. For instance, this is the case for many complex networks, such as power grids and biological networks. Complexity there arises from the interaction between network interconnectivity structure and the nonlinearities in the node and edge dynamics. And, this often leads to multistability. Linearization methods can provide insight near the system's attractors, but they poorly approximate the dynamics on the rest of the state space. 
In contrast, operator-based methods give access to the global characteristics of nonlinear systems. And, they do so in a linear setting and are therefore more well-suited, for instance, to characterize the attractor basin structure of multistable dynamical systems or the design control interventions. The Perron-Frobenius and Koopman operator are adjoint linear infinite-dimensional operators whose spectra can provide global information about the system. Their approximations using data-driven approaches make operator methods potentially applicable when there is no prior knowledge of the system.  

The Perron-Frobenius operator evolves densities on state space. It has been used extensively to access global behavior of nonlinear dynamical systems \cite{lasota1985probabilistic, tantet2018crisis}. There are several well developed approaches for obtaining its numerical approximations, such as Ulam's method that relies on the discretization of physical space to obtain an approximation of the Perron-Frobenius operator \cite{ulam1960collection}. Since the Koopman operator is adjoint to the Perron-Frobenius operator, numerical approximations of the Koopman operators can be obtained using these methods as well \cite{klus2015numerical}. 

The Koopman operator is an infinite dimensional linear operator that describes the evolution of \textit{observables} (functions of the state space) \cite{koopman1931hamiltonian, lasota1985probabilistic, klus2015numerical}. Its definition and properties in the context of dynamical systems are provided, for instance, in Ref. \cite{budivsic2012applied}, which \textcolor{black}{also summarizes its applicability to model reduction, coherency analysis, and ergodic theory}. Methods based on the Koopman operator decomposition have been proven useful for \textcolor{black}{applications such as model reduction and control of fluid flows \cite{bagheri2013koopman}, power system analysis \cite{susuki2016applied}, and extracting spatio-temporal patterns of neural data \cite{brunton2016extracting}}. 

Data driven methods to approximate the Koopman operator rely upon snapshot pairs of measurements of the state of the system at consecutive time steps. Reconstructing the operator from these snapshot pairs requires that a set of functions (called a dictionary of observables) is chosen. The first data driven method introduced, dynamic mode decomposition (DMD), implicitly uses linear monomials as a dictionary and thus is most applicable to systems where the Koopman eigenfunctions are well represented by this basis set \cite{schmid2010dynamic}. A more recent method called extended DMD (EDMD) introduced in Ref. \cite{williams2015data} can be more powerful than the standard DMD when applied to nonlinear systems as it allows the choice of more complicated sets of dictionary functions. Applying the EDMD is most computationally feasible if the number of the dictionary functions does not exceed the total number of the snapshot pairs used. That is not necessarily the case if a rich function dictionary (e.g., a dictionary of high order polynomials) is chosen. A modification of EDMD called kernel DMD introduced in Ref. \cite{williams2014kernel} addresses this issue by providing a way to efficiently calculate the Koopman operator approximation in a case when the number of dictionary functions exceeds the number of measurements.  Yet, the principled choice of an underlying dictionary that leads to an accurate approximation of the eigenspectrum corresponding to the leading Koopman modes using EDMD or kernel DMD remains an outstanding challenge. That problem is confronted, for instance, in Ref. \cite{li2017extended}, where an iterative EDMD dictionary learning method is presented. Although the optimal choice of dictionary functions is often unknown, there are some common choices that are known to produce accurate results for certain classes of systems \cite{williams2015data}.

Here we study nonlinear dynamical systems with discrete symmetries combining operator-based approaches and linear representation theory. Recently, related methods have been applied to dynamical systems with symmetries. On the one hand, Ref. \cite{mehta2006symmetry} addresses symmetries of the Perron-Frobenius operator in relation to the admissible symmetry properties of attractors. On the other, Ref. \cite{sharma2016correspondence} links the spatiotemporal symmetries of the Navier-Stokes equation to the spatial and temporal Koopman operator. Additionally, Ref. \cite{brunton2016discovering} noted that symmetry considerations play an important role in discovering governing equations. And, Ref. \cite{kaiser2018discovering} shows how conservation laws can be detected with Koopman operator approximations and then used to control Hamiltonian systems.

In contrast,  our focus is on dynamical systems with symmetries described by a finite group. We show how the properties of the associated Koopman operator spectrum can be linked to the properties of the spectrum of the finite dimensional approximations of the Koopman operator obtained from finite data. We further show how the analytic properties of the Koopman operator decomposition can inform the choice of dictionary functions that can be used in the Koopman operator approximations. This gives a practical way to reduce the dimensionality of the approximation problem. 

Our development builds as follows. \cref{sec:preliminaries} defines the Koopman operator, introduces approximation methods (EDMD and kernel DMD), and defines equivariant dynamical systems as well as useful concepts from group theory and representation theory. \cref{sec:Koopman} draws out the implications of dynamical system symmetries for the structure of the Koopman operator and its eigendecomposition. \cref{sec:EDMD} connects the properties of the Koopman operator and the structure of its EDMD approximation for symmetric systems. This then allows modifying the EDMD method to exploit the underlying symmetries, resulting in a block-diagonal Koopman operator approximation matrix. We also provide numerical examples, showing how using particular dictionary structures speed up the algorithm. Finally, the last section summarizes our findings and outlines directions for future work.

\section{Preliminaries} \label{sec:preliminaries}

\subsection{Koopman operator}

In this section, we provide some background in operator theoretic approaches to dynamical systems, in particular, the Koopman operator and its adjoint Perron-Frobenius operator. Since in this manuscript we address the approximations of the Koopman operator where the input is discrete time data, we focus on their definition in the discretized setting. The continuous time definitions are similar \cite{budivsic2012applied}. Our results regarding the degeneracy of Koopman operator eigenvalues and the properties of its corresponding eigenfunctions presented in \cref{sec:Koopman} hold in both discrete and continuous time settings. 

Suppose we are interested in studying continuous time autonomous dynamical systems defined as:
\begin{align}\label{eq:e1}
\dot{x}=g_c(x).
\end{align}
Here, $x\in \mathbb{R}^n$, $g_c:\mathbb{R}^n\rightarrow \mathbb{R}^n$.
Let $\Phi(x(t), \Delta t)$ be a flow map mapping the initial condition $x(t)$ to the solution at time $t+\Delta t$. It is defined in the following way:
\begin{align}
\Phi(x(t), \Delta t)=x(t)+\int\limits_{t}^{t+\Delta t} g_c(x(\tau)) d\tau.
\end{align}
The system can be discretized with a finite time step $\Delta t_{\text{step}}$, so that $x_{i+1}=\Phi(x_i,\Delta t_{\text{step}})$. We denote the function evolving the dynamics of this discretized system by $g$:
\begin{align}\label{eq:e2}
x_{i+1}=g(x_i).
\end{align}

The Koopman operator is a \textit{linear infinite dimensional} operator that evolves the functions (referred to as observables) of state space variables $f:\mathbb{R}^n\rightarrow \mathbb{C}$. The action of the Koopman operator $\mathcal{K}$ on an observable function $f$ for discrete time systems is defined as:
\begin{align}
(\mathcal{K} f)(x)=f(g(x)).
\end{align}
Since we consider data-driven Koopman operator approximation methods in this manuscript, the discrete time version of the definition is most applicable. 

In general, parts of the Koopman operator spectrum can be continuous \cite{korda2018data, budivsic2012applied}, for instance, that can be the case for chaotic systems. However, for the purposes of this manuscript, we focus on the case of a discrete spectrum. 

Pairs of eigenvalues $\lambda$ and eigenfunctions $\phi$ of the Koopman operator $\mathcal{K}$ are defined as:
\begin{align}
(\mathcal{K} \phi)(x)=\lambda\phi(x).
\end{align}
Of particular interest are the Koopman modes that can be used in model reduction and coherency estimation \cite{susuki2011nonlinear, rowley2017model}. The Koopman modes $v_i$ of the full state observable are defined by:
\begin{align}
    x=\sum\limits_i v_i\phi_i(x),
\end{align}
and are projections of the observable onto the span of the eigenfunctions of the Koopman operator $\mathcal{K}$. 

The other candidate for studying dynamical systems using an operator based approach is the Perron-Frobenius operator $\mathcal{P}$ defined as follows for deterministic dynamical systems:
\begin{align}
\int\limits_A \mathcal{P}\rho(x) dx = \int\limits_{g^{-1}(A)} \rho(x) dx.
\end{align}
Here, $\rho(x)$ is a density on state space, and $A\subseteq \mathbb{R}^n$ is a subset of the state space. The Perron-Frobenius operator is an adjoint to the Koopman operator \cite{lasota1985probabilistic}, so an approximation of one of them provides an approximation of the other \cite{ klus2015numerical}.

\subsection{Koopman operator approximation methods}

Extended dynamic mode decomposition (EDMD) introduced in Ref. \cite{williams2015data} is a way of approximating the Koopman operator for discretized systems that requires an explicit choice of a dictionary of functions referred to as \textit{observables}. How to optimally choose those functions remains an open problem for many systems, especially if the form of differential equations describing the governing dynamical system is not known in advance and only finite data on the behavior of the system is available. The method can be very accurate in capturing the dynamics of the system, but its accuracy depends strongly on the choice of an appropriate dictionary of observables. The method's convergence properties are studied in Ref. \cite{korda2018convergence}, and its relation to the Perron-Frobenius operator approximation methods is discussed in Ref. \cite{klus2015numerical}. Here, we summarize the EDMD and its relation to the Koopman operator. 

The first requirement for the method is a set pairs of consecutive snapshots $x=[x_1,x_2,...,x_M]$ and $y=[y_1,y_2,...,y_M]$, where the measurements $x_i$ and $y_i$ are performed with a small constant time interval $\Delta t$: $y_i=\Phi(x_i,\Delta t)$. Typically, the set of snapshots contains measurements from different trajectories in state space. We define a dictionary of linearly independent observables $\mathcal{D}=\{\psi_1,...,\psi_N\}$ and form vectors of observables $\Psi_x$ and $\Psi_y$. Here, $\Psi_x\in \mathbb{R}^{M\times N}$, where $N$ is the number of dictionary functions used in the approximation, and $M$ is the number of data snapshots. The elements of $\Psi_x$ are obtained from $(\Psi_x)_{ij}=\psi_j(x_i)$. We also use the notation $\Psi(x_m)=(\psi_1(x_m),...,\psi_N(x_m))$ for the dictionary functions evaluated at a particular point on the trajectory.

A finite dimensional approximation of the Koopman operator $\mathcal{K}$ that we denote as $K$ can be obtained using:
\begin{align}
K=\Psi_x^{+}\Psi_y.
\end{align}
Here, $\Psi_x^{+}$ denotes the pseudoinverse of $\Psi_x$. We focus on the case of the Moore-Penrose pseudoinverse for the rest of the manuscript \cite{penrose1955generalized}.

If the number of snapshots is much higher than the dimensionality of the function dictionary ($M\gg N$), it is more practical instead to define the square matrices $G$ and $A$ as shown below and obtain the approximation in the following way:
\begin{align}
\begin{gathered}
K=G^{+}A,\text{ where}\\ G=\sum\limits_m\Psi(x_m)^*\Psi(x_m),~~ A=\sum\limits_m\Psi(x_m)^*\Psi(y_m).
\end{gathered}
\end{align}
Here, $^*$ represents the complex conjugate transpose. If the only observables are the states of the system $x_1,x_2,...,x_n$, EDMD reduces to DMD \cite{williams2015data, klus2015numerical}. 

The eigendecomposition of $K$ provides the Koopman eigenvalues, eigenfunctions, and modes that allow an approximate linear representation of the underlying system dynamics. Let $\lambda_j$ and $u_j$ be the $j$\ts{th} eigenvalue and eigenvector of $K$. Then the corresponding Koopman eigenfunction can be approximated by:
\begin{align}
    \phi_j(x)=\Psi(x) u_j.
\end{align}
Let $b_i$ be the vectors defined by $g(x)_i=\Psi b_i$, where $g(x)_i=e^*_i x$ denotes the elements the full state observable discussed in Ref. \cite{williams2015data}, and $B=(b_1~~...~~b_n)$. The Koopman eigenmodes can then be obtained as:
\begin{align}
    v_i=(w_i^* B)^T.
\end{align}
Here, $w_i$ denotes the $i$\ts{th} left eigenvector of $K$.

A modification of EDMD named kernel DMD \cite{williams2014kernel} is better suited for systems with a low number of measurements and a high number of observables (e.g., the full state observable for fluid dynamical systems is very high dimensional, so defining a polynomial dictionary of the full state observable is very computationally expensive), i.e. $M\ll N$. The method relies on evaluating the kernel function: \begin{align}
    k(x_i,y_i)=\Psi(x_i)\Psi(y_i)^*.
 \end{align}
That allows efficient computation of $M\times M$ matrices $\hat{G}$, $\hat{A}$, and $\hat{K}$, where $M$ is the number of trajectory time steps. The eigendecomposition of $\hat{K}$ then can be used to obtain the approximations of the Koopman eigenvalues, eigenfunctions, and modes.

In the main body of the manuscript, we focus on the case where the number of measurements is relatively high for each degree of freedom ($M\gg N$), and obtain a way to reduce the dimensionality of the EDMD approximation of the Koopman operator for systems with symmetries in  \cref{sec:EDMD}. A similar modification of the kernel DMD is discussed in \cref{KDMD}.

\subsection{Point symmetries} 

In this section, we define the concepts useful to study the structure of the Koopman operator $\mathcal{K}$ and its approximations $K$ for systems with symmetries. Throughout this section and the rest of the manuscript, we use an example of a small network of Duffing oscillators to illustrate the definitions and algorithms. 

In this manuscript, we consider dynamical systems (as defined in \cref{eq:e1,eq:e2}) that respect point symmetries. These systems are called equivariant with respect to the symmetry group $\Gamma$. We define groups by their \textit{presentations} in a form $\langle S|R \rangle$, where $S$ is a set of generators of the group, and $R$ is a set of relations among these generators that define that group. Every element of the group can be written as a product of powers of some of these generators.

For instance, the cyclic group $Z_3$ is presented by $\langle r|r^n=1 \rangle$. An example of a realization of that group is a set of rotational symmetries of a regular $n$-gon.

To study dynamical systems with symmetries, we need to define the specific actions of the group on a vector space in addition to an abstract presentation of a group $\Gamma$. Let $X\subset\mathbb{R}^n$ be a vector space with elements $x\in X$. We denote the actions $\gamma_{\rho}$ on a vector space $X$ by $\gamma_{\rho} x$ if the set of these actions $\Gamma_{\rho}$ is isomorphic to $\Gamma$. A shorthand $\gamma_{\rho} x =\gamma x$ is sometimes used in the literature when the action corresponding to the subscript $\rho$ is clear from the context (for instance, it is defined by a permutation matrix of the same degree as the state space of the system), however, we use the $\gamma_\rho$ notation to avoid ambiguity, since the precise definition of group action in particular cases is important in this manuscript, as shown, for instance, in \cref{Ex:isotypic} and \cref{Ex:isotypic2}.

Finally, we define what it means for a dynamical system to be symmetric. Let $\dot{x} = g_c(x)$ be a continuous time system of differential equations. Here, $x\in \mathbb{R}^n$, and $g_c:\mathbb{R}^n \rightarrow \mathbb{R}^n$. The system is $\Gamma$-equivariant with respect to the actions of $\Gamma_\rho$ if for all $x\in\ X$ and $\gamma_{\rho} \in \Gamma_\rho$:
\begin{align}\label{equivariance}
g_c(\gamma_{\rho} x(t)) = \gamma_{\rho} g_c(x(t)).
\end{align}
As discussed in \cref{sec:preliminaries}, data comes in discretized form, so a discrete form of that definition is useful. For discrete time systems defined by $x_{i+1} = g(x_i)$, equivariance is defined in a similar manner:
\begin{align}
    g(\gamma_{\rho} x_i) = \gamma_{\rho} g(x_i).
\end{align}
We note that if a continuous time system is $\Gamma$-equivariant, so is its discretization. \textcolor{black}{Moreover, the set of trajectories of a $\gamma$-equivariant system in state space also respects the symmetries of the system}. For discretized systems, it means that if $\{x_0, x_1, ... x_n\}$ form a trajectory in state space, then $\{\gamma_{\rho}x_0, \gamma_{\rho}x_1, ... \gamma_{\rho}x_n\}$ form a trajectory as well.

\textcolor{black}{An important example of equivariant dynamical systems that many of the recent works have focused on (such as Refs.\cite{pecora2014cluster,hart2019topological,sorrentino2016complete, matheny2019exotic}) is a system of coupled identical oscillators. In that case, the set (or a subset) of actions under which the system is equivariant is defined by the set of permutational matrices $P$ that commute with the adjacency matrix of that oscillator network. In this case, the action of the group is linear, however, that does not always have to be the case.} 

We also need to define the action of the group in function space, where $f\in\mathcal{F}$ are functions $f:X\rightarrow\mathbb{C}$ as:
\begin{align}
(\gamma_{\rho}\circ f)(x) \equiv f(\gamma_{\rho}^{-1} x).
\end{align}
This definition will be useful since the Koopman operator acts on functions (i.e. observables).

Another concept useful for our work is a linear group \textit{representation} $T$, which is a mapping from group elements $\gamma\in\Gamma$ to the elements of the general linear group (a group of matrices of degree $n$ with the operation of matrix multiplication denoted by $\text{GL}(n,V)$) on a vector space $V$ (in this case, we are interested in $V=\mathbb{C}^n$). The characters of a group representation $T_i(\gamma)$ are defined as $\chi_i(\gamma)=\text{Tr} (T_i(\gamma))$. 

A representation is called irreducible if it has no nontrivial invariant subspaces (meaning that the representation matrices corresponding to the group elements can not be simultaneously non-trivially block diagonalized into the same block form). For each $\Gamma$ we can obtain all of its irreducible matrix representations. We denote their elements mapping $\gamma\in\Gamma$ to $p\times p$-dimensional matrices as $R_i(\gamma)$, where the index $i$ corresponds to the $i$\ts{th} irreducible representation. Irreducible representations are defined up to an isomorphism. For the purposes of this manuscript, it is useful to either make use of the unitary irreducible representations or the characters $\chi _i(\gamma)$ of the unitary irreducible representations. 

A vector space, e.g. the space of square integrable functions $\mathcal{F}$, can be uniquely decomposed into components that transform like the $i$\ts{th} irreducible representation of $\Gamma$ under the actions of $\Gamma_{\rho}$. These components are called \textit{isotypic components} \cite{golubitsky2003symmetry}. We denote these components by $\mathcal{F}_i$. An \textit{isotypic decomposition} of the square integrable function space with respect to $\Gamma_{\rho}$ is then defined as $\mathcal{F}=\bigoplus\limits_i \mathcal{F}_i$. We illustrate the construction of an isotypic decomposition using an example of a $Z_2$-equivariant system.

\begin{example} \label{Ex:isotypic} Symmetries of a single Duffing oscillator dynamics and isotypic components in
function space corresponding to the actions of its symmetry group.

The unforced Duffing oscillator equation has the form:
\begin{align}
\ddot{x}=-\sigma\dot{x}-x(\beta+\alpha^2 x).
\end{align}
We can rewrite the above equation as a system of differential equations to obtain:
\begin{align}
\begin{gathered}
\dot{x}=y,\\
\dot{y}=-\sigma y-x(\beta+\alpha^2 x).
\end{gathered}
\end{align}
Let $\bm{x}=\begin{pmatrix}
x\\y
\end{pmatrix}$, and let the dynamics be denoted by $\dot{\bm{x}}=g(\bm{x})$. Let $r_{s}=\begin{pmatrix}
-1&0\\0&-1
\end{pmatrix}=-I_{2\times 2}$ be the action on the state space that flips the signs of variables. The actions $r_{s}$ and $e_{s}=I_{2\times 2}$ form a group $\Gamma_{s}$ isomorphic to $\Gamma=Z_2=\langle r|r^2=e\rangle$. Let $\gamma_{s}\in\Gamma_{s}$, then:
\begin{align}
\gamma_{s} g(\bm{x}) = g(\gamma_{s} \bm{x}).
\end{align}
Thus, the Duffing oscillator system is $Z_2$-equivariant with respect to the actions $\gamma_{s}$. 

We now illustrate the isotypic component decomposition of $Z_2$ in function space. $Z_2$ has two one-dimensional irreducible representations: the \textit{trivial} representation defined by $R_{\text{tr}}(r)=1$ and the \textit{sign} representation defined by $R_{\text{sign}}(r)=-1$. 
Then the space of square integrable functions $\mathcal{F}$ can be decomposed into $\mathcal{F}=\mathcal{F}_{\text{tr}}\oplus\mathcal{F}_{\text{sign}}$, where $\mathcal{F}_{\text{tr}}=\{f:r_{s}\circ f=f(-x,-y)=f(x,y)\}$ and $\mathcal{F}_{\text{sign}}=\{f:r_{s}\circ f=f(-x,-y)=-f(x,y)\}$.
In this case, the sets of functions $\mathcal{F}_{\text{tr}}$ and $\mathcal{F}_{sign}$ consisting of even and odd functions respectively transform like the trivial and sign irreducible representations with respect to sign flip as a group generator action.
\end{example}
We now extend the example to a network of Duffing oscillators and explore additional permutation symmetries. 

\begin{example} \label{Ex:isotypic2}
We now consider the dynamics of a network of Duffing oscillators. Suppose the coupling is linear in $x$ with a coupling coefficient assigned to every edge $\eta_{ij}$. Then for each node $i$ in the network, we have the following dynamics:
\begin{align}
\begin{gathered}
\dot{x_i}=y_i,\\
\dot{y_i}=-\sigma\dot{y_i}-x_i(\beta+\alpha^2 x_i)+\sum\limits_{ij}\eta_{ij}(x_i-x_j).
\end{gathered}
\end{align}
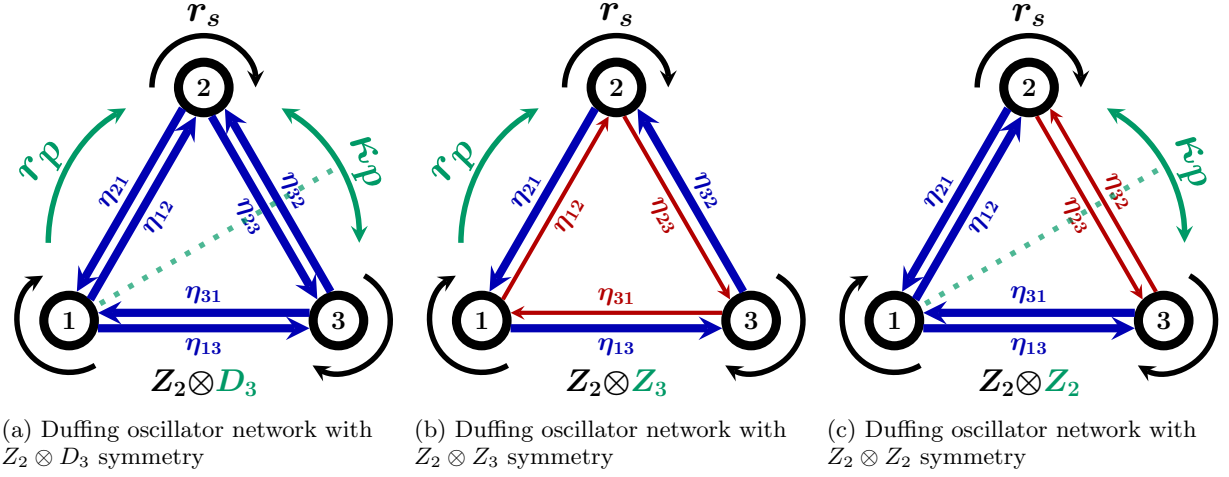
\begin{figure*}
%\label{fig:Duff Net}
\begin{subfigure}[b]{0.3\textwidth}
\centering
        \resizebox{\linewidth}{!}{
\begin{tikzpicture}
\tikzset{edge/.style = {->,> = stealth, line width = .5}}
\tikzstyle{every node}=[font=\Large]
\tikzset{vertex/.style = {shape=circle,fill=white,text=black,line width=5,draw,minimum size=3.em}}
\node[vertex] (B) at (210:3) {\textbf{1}};
\node[vertex] (A) at (90:3) {\textbf{2}};
\node[vertex] (C) at (330:3) {\textbf{3}};
\draw[edge, line width = 4.5, blue!70!black] (A.225) -- node[rotate=60,above] {$\bm{\eta_{21}}$} (B.75);
\draw[edge, line width = 4.5, blue!70!black] (B.45) -- node[rotate=60,below] {$\bm{\eta_{12}}$} (A.255);
\draw[edge, line width = 4.5, blue!70!black] (C.105) -- node[rotate=300,above] {$\bm{\eta_{32}}$} (A.315);
\draw[edge, line width = 4.5, blue!70!black] (A.285) -- node[rotate=300,below] {$\bm{\eta_{23}}$} (C.135);
\draw[edge, line width = 4.5, blue!70!black] (B.345) -- node[rotate=0,below] {$\bm{\eta_{13}}$} (C.195);
\draw[edge, line width = 4.5, blue!70!black] (C.165) -- node[rotate=0,above] {$\bm{\eta_{31}}$} (B.15);

\node[font=\LARGE] at (0,-2.7) {$\bm{\textcolor{black}{Z_2\otimes} \textcolor{green!60!blue}{D_3}}$};

\node[] (D) at (0,0) {};
\draw[edge, line width = 3., <->, opacity=1.]  (D) pic[green!60!blue]{carc=60:0:3.cm};
\node [rotate=-60, green!60!blue] at ({3.1},{1.55}) {\fontsize{30}{30}\selectfont \textbf{$\bm{{\kappa_{p}}}$}};

\draw[edge, line width = 3., ->, opacity=1.]  (D) pic[green!60!blue]{carc=180:120:3.cm};
\node [rotate=60, green!60!blue] at ({-3.1},{1.55}) {\fontsize{30}{30}\selectfont \textbf{$\bm{{r_{p}}}$}};

\draw[edge, line width = 3.]  (A) pic[black] {carc=180:0:1.cm};
\node [rotate=0, black] at ({0},{4.4}) {\fontsize{20}{20}\selectfont \textbf{$\bm{{r_{s}}}$}};
\draw[edge, line width = 3.]  (B) pic[black] {carc=-60:-240:1.cm};
\draw[edge, line width = 3.]  (C) pic[black] {carc=60:-120:1.cm};

\draw [dashed, loosely dashed, line width = 3., green!60!blue, opacity = 0.7] (-2,-1.175) -- (2.6,1.475);
\end{tikzpicture}
}
\caption{Duffing oscillator network with $Z_2\otimes D_3$ symmetry}
\label{subfig:Z2_D3}
\end{subfigure}
%%%%%%%%%%%%%%%%%%%%%%%%%%%%%%%%%%%%%%%%%%%%%%%%%%%%%%%
     \begin{subfigure}[b]{0.3\textwidth}
     \centering
        \resizebox{\linewidth}{!}{
\begin{tikzpicture}
\tikzset{edge/.style = {->,> = stealth, line width = .5}}
\tikzstyle{every node}=[font=\Large]
\tikzset{vertex/.style = {shape=circle,fill=white,text=black,line width=5,draw,minimum size=3.em}}
\tikzset{vertex/.style = {shape=circle,fill=white,text=black,line width=5,draw,minimum size=3.em}}
\node[vertex] (B) at (210:3) {\textbf{1}};
\node[vertex] (A) at (90:3) {\textbf{2}};
\node[vertex] (C) at (330:3) {\textbf{3}};
\draw[edge, line width = 4.5, blue!70!black] (A.225) -- node[rotate=60,above] {$\bm{\eta_{21}}$} (B.75);
\draw[edge, line width = 2.3, red!70!black] (B.45) -- node[rotate=60,below] {$\bm{\eta_{12}}$} (A.255);
\draw[edge, line width = 4.5, blue!70!black] (C.105) -- node[rotate=300,above] {$\bm{\eta_{32}}$} (A.315);
\draw[edge, line width = 2.3, red!70!black] (A.285) -- node[rotate=300,below] {$\bm{\eta_{23}}$} (C.135);
\draw[edge, line width = 4.5, blue!70!black] (B.345) -- node[rotate=0,below] {$\bm{\eta_{13}}$} (C.195);
\draw[edge, line width = 2.3, red!70!black] (C.165) -- node[rotate=0,above] {$\bm{\eta_{31}}$} (B.15);

\node[font=\LARGE] at (0,-2.7) {$\bm{\textcolor{black}{Z_2\otimes} \textcolor{green!60!blue}{Z_3}}$};

\draw[edge, line width = 3.]  (A) pic[black] {carc=180:0:1.cm};
\node [rotate=0, black] at ({0},{4.4}) {\fontsize{20}{20}\selectfont \textbf{$\bm{{r_{s}}}$}};
\draw[edge, line width = 3.]  (B) pic[black] {carc=-60:-240:1.cm};
\draw[edge, line width = 3.]  (C) pic[black] {carc=60:-120:1.cm};

\node[] (D) at (0,0) {};

\draw[edge, line width = 3., ->, opacity=1.]  (D) pic[green!60!blue]{carc=180:120:3.cm};
\node [rotate=60, green!60!blue] at ({-3.1},{1.55}) {\fontsize{30}{30}\selectfont \textbf{$\bm{{r_{p}}}$}};

\end{tikzpicture}
}
\caption{Duffing oscillator network with $Z_2\otimes Z_3$ symmetry}
\label{subfig:Z2_Z3}
\end{subfigure}
\begin{subfigure}[b]{0.3\textwidth}
\centering
        \resizebox{\linewidth}{!}{
\begin{tikzpicture}
\tikzset{edge/.style = {->,> = stealth, line width = .5}}
\tikzstyle{every node}=[font=\Large]
\tikzset{vertex/.style = {shape=circle,fill=white,text=black,line width=5,draw,minimum size=3.em}}
\node[vertex] (B) at (210:3) {\textbf{1}};
\node[vertex] (A) at (90:3) {\textbf{2}};
\node[vertex] (C) at (330:3) {\textbf{3}};
\draw[edge, line width = 4.5, blue!70!black] (A.225) -- node[rotate=60,above] {$\bm{\eta_{21}}$} (B.75);
\draw[edge, line width = 4.5, blue!70!black] (B.45) -- node[rotate=60,below] {$\bm{\eta_{12}}$} (A.255);
\draw[edge, line width = 2.3, red!70!black] (C.105) -- node[rotate=300,above] {$\bm{\eta_{32}}$} (A.315);
\draw[edge, line width = 2.3, red!70!black] (A.285) -- node[rotate=300,below] {$\bm{\eta_{23}}$} (C.135);
\draw[edge, line width = 4.5, blue!70!black] (B.345) -- node[rotate=0,below] {$\bm{\eta_{13}}$} (C.195);
\draw[edge, line width = 4.5, blue!70!black] (C.165) -- node[rotate=0,above] {$\bm{\eta_{31}}$} (B.15);

\node[font=\LARGE] at (0,-2.7) {$\bm{\textcolor{black}{Z_2\otimes} \textcolor{green!60!blue}{Z_2}}$};

\node[] (D) at (0,0) {};
\draw[edge, line width = 3., <->, opacity=1.]  (D) pic[green!60!blue]{carc=60:0:3.cm};
\node [rotate=-60, green!60!blue] at ({3.1},{1.55}) {\fontsize{30}{30}\selectfont \textbf{$\bm{{\kappa_{p}}}$}};

\draw[edge, line width = 3.]  (A) pic[black] {carc=180:0:1.cm};
\node [rotate=0, black] at ({0},{4.4}) {\fontsize{20}{20}\selectfont \textbf{$\bm{{r_{s}}}$}};
\draw[edge, line width = 3.]  (B) pic[black] {carc=-60:-240:1.cm};
\draw[edge, line width = 3.]  (C) pic[black] {carc=60:-120:1.cm};

\draw [dashed, loosely dashed, line width = 3., green!60!blue, opacity = 0.7] (-2,-1.175) -- (2.6,1.475);
\end{tikzpicture}
}
\caption{Duffing oscillator network with $Z_2\otimes Z_2$ symmetry}
\label{subfig:Z2_Z2}
\end{subfigure}
\caption{Possible symmetries of a network of three identical Duffing oscillators depending on the coupling strength between the oscillators. Different coupling strengths are shown in blue and red. Green arrows correspond to permutational symmetries arising from physical coupling, black arrows correspond to the symmetries of nodal dynamics. }
\label{fig:Duff Net}
\end{figure*}

This general coupling scheme is used to describe many systems in literature \cite{pecora2014cluster, sorrentino2016complete}. 

We now consider the case of a 3-node network. Depending on what the coupling terms are, the system may be $\Gamma$-equivariant with respect to different symmetry groups that act by permuting node indexes. Some examples are:
\begin{enumerate}[label=$\alph*)$]
\item{If all coupling strengths $\eta_{ij}$ are equal, the network has $D_3$ symmetry. This case is shown on \cref{subfig:Z2_D3}. Let the state of the system be defined by $x=\left( x_1~~y_1~~x_2~~y_2~~x_3~~y_3\right)^T$.
Then, the symmetry group is presented by $D_3=\langle r,\kappa | r^3=\kappa^2=e, \kappa r \kappa = r^{-1} \rangle$ and generated by the actions $r_p=\begin{pmatrix}
0&1&0\\
0&0&1\\
1&0&0
\end{pmatrix} \otimes I_{2\times 2}$ and $\kappa_p=\begin{pmatrix}
1&0&0\\
0&0&1\\
0&1&0
\end{pmatrix} \otimes I_{2\times 2}$.}
\item{If the coupling strengths obey the conditions $\eta_{ij}\neq \eta_{ji}$ and $\eta_{ij}=\eta_{jk}$, the network has $Z_3$ symmetry. This case is shown on \cref{subfig:Z2_Z3}.
The symmetry group is presented by $Z_3=\langle r| r^3=e \rangle$ and generated by the action $r_p$ defined above.}
\item{If the coupling strengths obey the conditions $\eta_{12} = \eta_{21}$ and no other equalities hold, the network has $Z_2$ symmetry. This case is shown on \cref{subfig:Z2_Z2}.
The symmetry group is presented by $Z_2=\langle \kappa| \kappa^2=e \rangle$ and generated by the action $\kappa_p$ defined above.}
\end{enumerate}
Even though in case (c) the permutation symmetry is isomorphic to the same group as the sign flip symmetry in \cref{Ex:isotypic}, the isotypic components in function space $\mathcal{F}$ induced by the group action are different. $Z_2$ has two one-dimensional irreducible representations: \textit{trivial} representation $R_1(r)=R_{tr}(r)=1$ and \textit{sign} representation $R_2(r)=R_{sign}(r)=-1$. Let $\bm{x_i}=\begin{pmatrix}
x_i\\y_i
\end{pmatrix}$. The isotypic components are defined by $\mathcal{F}_{tr}=\{f:r_{p}\circ f=f(\bm{x}_1,\bm{x}_3,\bm{x}_2)=f(\bm{x}_1,\bm{x}_2,\bm{x}_3)\}$ and $\mathcal{F}_{sign}=\{f:r_{p}\circ f=f(\bm{x}_1,\bm{x}_3,\bm{x}_2)=-f(\bm{x}_1,\bm{x}_2,\bm{x}_3)\}$.

Additionally, each node still has $Z_2$ symmetry with respect to the action $r_s$ which is not broken since the coupling function is odd. That symmetry is also depicted in \cref{fig:Duff Net}. We can find the isotypic components of the whole symmetry group as $\mathcal{F}=\bigoplus\limits_{i,j} \left( \mathcal{F}_{i,s}\otimes \mathcal{F}_{j,p}\right)$.
\end{example}

Any function can be rewritten as a sum of projections into different isotypic components. The procedure is outlined in the following section.

\section{Properties of the Koopman operator for systems with symmetries}\label{sec:Koopman}

In this section, we consider the structure of the eigenspace of the Koopman operator of $\Gamma$-equivariant systems. We show how to obtain a particular eigenbasis of the system corresponding to the isotypic decomposition in function space and demonstrate that the isotypic decomposition induces a block diagonal structure on the matrix representation of $\mathcal{K}$.

\begin{theorem} \label{Comm} For a $\Gamma$-equivariant dynamical system $\dot{x} = g(x)$ and an arbitrary function $f$, the Koopman operator commutes with the actions of the elements of $\Gamma$:
\begin{align}
\gamma_{\rho} \circ (\mathcal{K} f) (x)  = \mathcal{K} (\gamma_{\rho} \circ f)(x).
\end{align}
\textbf{Proof:}\\
The commutativity follows from the definitions of the Koopman operator and the definition of the action of the group in state space and function space.
\begin{align}
\begin{gathered}
\gamma_{\rho} \circ (\mathcal{K} f) (x) = \gamma_{\rho}\circ f(g(x))= f(\gamma_{\rho}^{-1}g(x))\\=
f(g(\gamma_{\rho}^{-1}x))=\mathcal{K}f(\gamma_{\rho}^{-1} x)=\mathcal{K} (\gamma_{\rho} \circ f)(x)
\end{gathered}
\end{align}
This result is similar to Theorem 3.1 in \cite{mehta2006symmetry}, where it is shown that the action of Perron-Frobenius operator commutes with the action of the symmetry group $\Gamma$ for $\Gamma$-equivariant systems.
\end{theorem}
\begin{corollary} \label{Eig} The space of eigenfunctions of the Koopman operator $\mathcal{K}$ with eigenvalue $\lambda$ for a $\Gamma$-equivariant system is $\Gamma$-invariant.\\
\textbf{Proof:}\\
Let $S_{\lambda}$ be the set of eigenfunctions of $\mathcal{K}$ with eigenvalue $\lambda$. Let $\phi\in S_{\lambda}$. Then, using the commutativity of $\mathcal{K}$ and $\Gamma_{\rho}$, we can show that $\forall\gamma_{\rho}\in\Gamma_{\rho}$:
\begin{align}
\mathcal{K} (\gamma_{\rho} \circ \phi(x)) = \gamma_{\rho} \circ (\mathcal{K} \phi)(x) = \lambda \gamma_{\rho} \circ \phi(x).
\end{align}
Thus, $\phi_{\gamma,\rho}\in S_{\lambda}$, where $\phi_{\gamma,\rho}$ is also an eigenfunction with an eigenvalue $\lambda$ defined as $\phi_{\gamma,\rho}(x)=\gamma_{\rho} \circ \phi(x)$.
\end{corollary}

We now consider a particular form of the eigenbasis of the Koopman operator that induces block diagonal structure of the matrix representation of the action of the Koopman operator $\mathcal{K}$. The result quoted below is useful for that purpose. 

\begin{theorem}(Theorem 3.5 in Chapter XII of  Ref.\cite{golubitsky2012singularities}).\\
Let $\Gamma$ be a compact Lie group acting on the vector space $V$ decomposed into isotypic components $V=W_1\bigoplus...\bigoplus W_t$. Let $A:V\rightarrow V$ be a linear mapping commuting with $\Gamma$. Then $A(W_k)\subset W_k$. 
\end{theorem}

This result is applicable to finite symmetry groups. Isotypic components of $\mathcal{F}$ with respect to $\Gamma$ induce block diagonal structure of the matrix representation of the Koopman operator. Since $\mathcal{K}$ and $\Gamma$ commute, $\mathcal{K}(\mathcal{F}_k)\subset 
\mathcal{F}_k$. This block structure can be exploited in finding the Koopman operator approximations, as we show in the next section. Thus, we need to be able to obtain an isotypic component basis from an arbitrary function dictionary. This is a well defined procedure \cite{cornwell1997group}, outlined below. Functions obtained via isotypic decomposition are useful to perform calculations in many areas of physics, for instance, they can simplify finding approximate solutions to Schrodinger equation, or in studying crystallographic point groups \cite{cornwell1997group, stiefel2012group}. The construction is also widely applied to dynamical systems, for instance, to study states and their stability using equivariant bifurcation theory.

Suppose we start from an arbitrary basis function dictionary $\mathcal{D}_{\psi}=\{\psi_i\}$. Each of those functions can be expanded in the isotypic component basis with at least one nonzero coefficient $\alpha^p_{mn}$:
\begin{align}
\psi=\sum\limits_{p}\sum\limits_{m,n=1}^{d_p} \alpha^p_{mn}\xi^p_{mn}.
\end{align}
Here, $\xi^p_{mn}$ is a basis function in the $p$\ts{th} isotypic component of $\mathcal{F}$, and $d_p$ is the dimension of that isotypic component. Alternatively, it can be thought of as a sum over all inequivalent (non-isomorphic) irreducible representations of $\Gamma$, where $\xi_{mn}^p$ transforms as $(m,n)$\ts{th} element of $p$\ts{th} irreducible representation of $\Gamma$ \cite{cornwell1997group}. We define a projection operator and form a new function basis consisting of functions $\{\xi_{mn}^p\}$ as outlined below. 

The projection operator is defined as:
\begin{align} \label{eq:projection operator1}
\mathcal{P}^p_{mn}=\dfrac{d_p}{|\Gamma|}\sum\limits_{\gamma\in\Gamma} [R_p(\gamma)]^*_{mn} \gamma_{\rho}.
\end{align}
Here, $[R_p(\gamma)]_{mn}$ denotes the element in $n$\ts{th} row and $m$\ts{th} column of the $i$ th unitary irreducible representation of $\gamma\in\Gamma$, and $\gamma_{\rho}$ is the group action. We can form an orthonormal basis $\mathcal{D}_{\xi}=\{\xi_i\}$ using that projection operator as follows:
\begin{align} \label{eq:basis}
\xi^p_{mn}(x)=\dfrac{1}{c_{np}}\mathcal{P}^p_{mn}\circ\psi(x).
\end{align} 
Here, $c_{np}=\langle \mathcal{P}^p_{nn} \psi, \mathcal{P}^p_{nn} \psi \rangle^{1/2}$, where $\langle , \rangle$ denotes the inner product\textcolor{black}{, which can be omitted for our purposes since the scaling of basis functions does not affect the EDMD results.}

Equivalently, due to orthogonality of characters of irreducible representations, the projection operator can be obtained using the following expression:
\begin{align} \label{eq:projection operator}
\mathcal{P}^p=\dfrac{d_p}{|\Gamma|}\sum\limits_{\gamma\in\Gamma} \chi_p(\gamma)^*\gamma_{\rho}.
\end{align}
Here, $\chi_p(\gamma)$ is a character of the $p$\ts{th} irreducible representation of $\Gamma$.
If this formula is used, each irreducible representation of degree $d_p$ provides a basis function, and $d_p^2-1$ other basis functions can be formed using the Gram-Schmidt orthogonalization process \cite{cornwell1997group, stiefel2012group, pecora2014cluster}. 

Once an isotypic component basis is obtained, the action of the Koopman operator on function space can be presented in the form of a block diagonal matrix. Each irreducible unitary representation of dimension $d_p$ in this case corresponds to a number $d_p$ of $d_p\times d_p$ sized blocks in that matrix $\mathcal{K}$. Similar analysis works for the Koopman operator approximation $K$. \textcolor{black}{The reason why this additional decomposition works can be found in  \cref{app:iso}.} 

\section{Implications for EDMD} \label{sec:EDMD}

\subsection{Constructing a basis for systems with known symmetries}

In this section we show that the approximation of $K$ obtained using EDMD can be reduced to the block diagonal structure similar to $\mathcal{K}$ under certain assumption on the data. We provide some examples of constructing an isotypic component basis from a given function dictionary. We highlight that the basis depends on both the structure of $\Gamma$ and the definition of its actions $\Gamma_{\rho}$.

First, we establish that the Koopman operator approximation $K$ commutes with the actions $\gamma_{\rho}$ of $\Gamma$ if the data used in the calculation respects the symmetry, meaning the set of pairs of data points satisfies the condition:
\begin{align}\label{eq:data symm}
\{(\gamma_{\rho} x_i, \gamma_{\rho} y_i)\}=\{(x_i,y_i)\}.
\end{align}
In other words, the set of trajectories is closed under the action of the symmetry group of the underlying dynamical system.

In order to perform further simplifications, we pick a particular order of group elements $\{\gamma_1,..., \gamma_{|\Gamma|}\}$ and create the vectors $\Psi_x$ (and analogously $\Psi_y$) according to that ordering:
\begin{align} \label{eq: order}
\begin{gathered}
    \Psi(x)=\begin{pmatrix}
    \gamma_1 \circ \Psi(x)~~...~~
    \gamma_{|\Gamma|} \circ \Psi(x)
    \end{pmatrix},\\
    \Psi_x=\begin{pmatrix}
    \Psi_1(x_1)&...&\Psi_{N/ |\Gamma|}(x_1)\\
    \vdots&\ddots&\vdots\\
    \Psi_1(x_M)&...&\Psi_{N/ |\Gamma|}(x_M)
    \end{pmatrix}.
\end{gathered}
\end{align}

Given the ordering of the group elements, we can also construct the permutation representation of the group such that: 
\begin{align} \label{eq:permutation}
    P_{\gamma_k} (\gamma_1,...,\gamma_{|\Gamma|})^T = (\gamma_k \gamma_1,...,\gamma_k \gamma_{|\Gamma|})^T.
\end{align}
By Cayley's theorem, such permutations form a group isomorphic to $\Gamma$. Determining the actions $P_{\gamma_k}$ of the group generators is sufficient to find the actions of all group elements. Let $\bm{P}_{\gamma_k}=P_{\gamma_k}\otimes I_{n\times n}$. We note that $(\bm{P}_{\gamma_k})^*=(\bm{P}_{\gamma_k})^{-1}$. It can be shown that:
\begin{align}\label{eq:comm}
    \bm{P}_{\gamma_k} G = G \bm{P}_{\gamma_k},~~\bm{P}_{\gamma_k} A = A \bm{P}_{\gamma_k}.
\end{align}
By definition, $A=\Psi_x^*\Psi_y$. We note that:
\begin{align}
    (\Psi_x\bm{P}_{\gamma_k})^* \Psi_y\bm{P}_{\gamma_k}= \Psi_x^*\Psi_y.
\end{align}
That is the case for symmetric trajectories satisfying \cref{eq:data symm}:
\begin{align}
\begin{gathered}
\left((\Psi_x\bm{P}_{\gamma_k})^* \Psi_y\bm{P}_{\gamma_k}\right)_{ij}=\sum\limits_m\psi_i^*(\gamma_k^{-1} x_m)\psi_j(\gamma_k^{-1} x_m)\\=
\sum\limits_m\psi_i^*(x_m)\psi_j(x_m)=
\left(\Psi_x^* \Psi_y\right)_{ij}.
\end{gathered}
\end{align}
Thus, $A$ and $G$ commute with the action of the symmetry group.

If $G$ is invertible and $G$ commutes with $\gamma_{\rho}$, $G^{-1}$ commutes with $\gamma_{\rho}$ as well. Then:
\begin{align}\label{eq: comm}
\begin{gathered}
\bm{P}_{\gamma_i}K = \bm{P}_{\gamma_i} G^{-1}A= G^{-1}A \bm{P}_{\gamma_i} = K \bm{P}_{\gamma_i}.
\end{gathered}
\end{align}

If $G$ is not invertible, the commutativity result still holds for $G^{+}$. $G$ is a normal matrix since it satisfies $GG^*=G^*G$. In \cref{app: K commutes}, we show that if $G$ is normal, $GG^{+}=G^{+}G$, so $G$ commutes with its Moore-Penrose pseudoinverse, and therefore the actions of $K$ and $\Gamma$ commute.

Since $K$ commutes with the actions of $\Gamma$, $K\mathcal{F}_i\subset \mathcal{F}_i$. This shows that $K$ can be block-diagonalized in the same way as $\mathcal{K}$.

Suppose we start from a dictionary of observables. Since that dictionary is not necessarily an isotypic component dictionary corresponding to $\Gamma$ and its action, in order to obtain a block diagonal matrix $K$, the dictionary needs to be modified using the procedure outlined in \cref{sec:Koopman}. In the example below, we show explicitly how to perform this transformation into the isotypic component basis. 

In order for the basis to faithfully represent the symmetries of the system we require that:
\begin{itemize}
\item
The dictionary is closed under the action of the symmetries of the system:
\begin{align}
    \text{If } \psi\in\mathcal{D},~\gamma_{\rho}\psi\in \text{span}(\mathcal{D})
\end{align} 

\item
Each isotypic component is present after the isotypic component decomposition of the original function basis:
\begin{align}
    \forall m,p~~ \exists \psi\in\mathcal{D}, \text{ s.t. } \mathcal{P}^p_{mn}\psi\neq 0
\end{align}
\end{itemize}

For instance, using a monomial basis for a $D_3$ equivariant system does not satisfy the second requirement.

\textcolor{black}{If these requirements are satisfied, the change of basis does not affect the result obtained by applying the EDMD algorithm as shown in \cref{app: basis}.}  Additionally, we note that the eigenvalues of $K$ do not typically have the same degeneracy properties as the eigenvalues of $\mathcal{K}$, but the symmetries of the underlying dynamical system are preserved in trajectory reconstructions.

\begin{example}\label{ex:iso Duff}
Constructing an isotypic component basis for a network of Duffing oscillators from a given basis. 

We consider the case of a system of Duffing oscillators with identical coupling as depicted in \cref{subfig:Z2_D3}. In that case, the system has $Z_2\otimes D_3$ symmetry. Suppose we want to construct an isotypic component basis from a given function dictionary $\mathcal{D}$. As an example, we use an initial dictionary $\mathcal{D}_{\psi}$ of $n$ mesh-free radial basis functions. The radial basis function centers can be obtained by either k-means clustering of the data or sampling from a predetermined distribution. Each function can be presented in a form \textcolor{black}{$\psi(c,x)=r_{c} log (r_{c})$, where $c$ is a 6-dimensional radial basis function center, and $r_{c,x}=||x-c||^{1/2}$}. In order to preserve the symmetries of the system, we need to have dictionary elements corresponding to acting on the basis functions by each $\gamma_{\rho}\in\Gamma_{\rho}$. Due to the form of these functions, $\gamma_{\rho}\circ\psi(c,x)=\psi(\gamma_{\rho}^{-1}c,x)$.
 
Since $\Gamma=Z_2\otimes D_3$ is a direct product of two groups, we can write the projection operator in the following form:

\begin{align}
\mathcal{P}^p_{mn}=\dfrac{d_p}{|\Gamma|} \sum\limits_{\gamma\in\Gamma} |R_p(\gamma_i,\gamma_j)|^*_{mn}(\gamma_{i},\gamma_j).
\end{align} 

We first consider the irreducible representations of $D_3$. 
\begin{itemize}
\item{}Trivial representation $R_{tr}$: $R_{tr}(r)=1$, $R_{tr}(\kappa)=1$
\item{}Sign representation $R_{sign}$: $R_{sign}(r)=1$, $R_{sign}(\kappa)=-1$
\item{}Standard representation $R_{st}$: 
$R_{st}(r)=\begin{pmatrix}
\omega&0\\0&\omega^2
\end{pmatrix},~~
R_{st}(\kappa)=
\begin{pmatrix}
0&1\\1&0
\end{pmatrix}$.
Here, $\omega=e^{2\pi i/3}$.
\end{itemize}
Suppose we form a vector of basis functions in $\mathcal{D}_{\psi}$, $\Psi=(\psi_{1,1}\psi_{2,1},...,\psi_{n,1},...,\psi_{1,|\Gamma|},\psi_{2,|\Gamma|},...,\psi_{n,|\Gamma|})^T$, where the first index corresponds to acting on $\psi_{1,i}$ by the $j^{th}$ element of $\Gamma_{\rho}$. 
Using the equations above, we obtain transformation matrices that we will use to get the isotypic component basis: 
\begin{align}
T_{D_3}=\begin{pmatrix}
1&1&1&1&1&1\\
1&1&1&-1&-1&-1\\
1&\omega&\omega^2&0&0&0\\
0&0&0&1&\omega^2&\omega\\
1&\omega^2&\omega&0&0&0\\
0&0&0&1&\omega&\omega^2\\
\end{pmatrix},~~T_{Z_2}=
\begin{pmatrix}
1&1\\1&-1
\end{pmatrix}.
\end{align}
The isotypic component basis then can be obtained by modifying a set of functions in $\mathcal{D}_{\psi}$:
\begin{align}
\Xi=T \Psi,~~T=T_{Z_2}\otimes T_{D_3}\otimes I_{n\times n}.
\end{align}

If we use $\Xi$ as a basis, we obtain $K$ decomposed into 8 blocks, each corresponding to an irreducible representation of $Z_2\otimes D_3$. 

We implement the EDMD algorithm to obtain the approximation of $\mathcal{K}$. Here, the data comes from 500 initial trajectories of length 10 that were then reflected and rotated so that the data respects the symmetries. The parameter values of $\alpha=1.$, $\beta=-1.$, $\delta=.5$, and $\eta=1.$ were used. We plot the approximation matrix $K$ on \cref{fig: block}. In this case, a dictionary of \textcolor{black}{120 radial basis functions} was used. \cref{subfig: standard} illustrates the Koopman operator approximation matrix $K$ calculated using an initial dictionary $D_{\psi}$ and requires performing matrix operations on the full $120\times 120$ matrix. \cref{subfig: symm} shows $K$ obtained from the symmetry adapted basis functions. The order of calculations can be reduced significantly since it is only necessary to perform matrix operations on blocks. $K_{symm}$ has \textcolor{black}{$4$, $10\times 10$ and $2$, $20\times 20$} unique blocks.
\begin{figure}
\begin{subfigure}[b]{.23\textwidth} 
\includegraphics[width=0.9\textwidth]{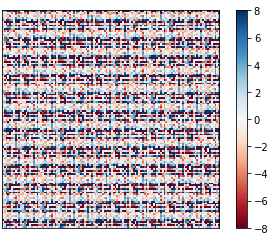}
\caption{$K$ for a standard dictionary of\\ observables}
\label{subfig: standard}
\end{subfigure}
\begin{subfigure}[b]{0.23\textwidth}
\includegraphics[width=0.9\textwidth]{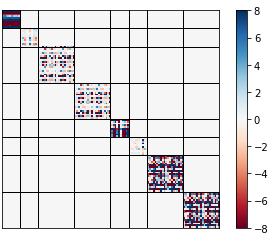}
\caption{$K$ for a symmetry adapted dictionary of observables}
\label{subfig: symm}
\end{subfigure}
\caption{Structure of $K$ for different choices of dictionary functions}
\label{fig: block}
\end{figure} 
\end{example} 

As shown in examples above, we can construct a basis that diagonalizes the Koopman operator matrix approximation $K$ from the elements of any arbitrary basis. Since the off block-diagonal elements of the matrix are a priori known to be zero, we do not need to compute these elements explicitly. This suggests that for systems with symmetries it is more efficient to perform the EDMD algorithm for isotypic decomposition blocks. We denote the number of conjugacy classes or irreducible representations of $\Gamma$ by $r_{\Gamma}$. In that case, instead of performing $O((mr_{\Gamma})^{\alpha})$ operations of matrix inversion, multiplication, and eigendecomposition, it is sufficient to perform these operations for each of the $r_{\Gamma}$ blocks, with operations being $O(m^{\alpha})$. Here $2<\alpha<3$ \cite{knuth1997art}. Even though the algorithmic complexity only differs by a factor that scales with the size of the group that is fixed for any given system, in practice, the computation is more efficient when EDMD specific to $\Gamma$-equivariant systems is used. We also note that each $d_p$ dimensional irreducible representation results in $d_p$ of $d_p\times d_p$ dimensional blocks in $K_{\text{symm}}$ known to be equal a priori, which further simplifies the calculation. Moreover, in the case of networks of high dimensionality, it allows parallel eigendecomposition computation for isotypic component blocks.  \cref{ttable} summarizes the modified EDMD algorithm for $\Gamma$-equivariant systems and highlights that the order of computations can be lowered.

\begin{table*} 
\fbox{%
    \parbox{1.\textwidth}{
\begin{minipage}{0.25\linewidth} 
\centering
{\textbf{Standard EDMD}} 
\end{minipage}
\begin{minipage}{0.73\linewidth}
\centering
{\textbf{EDMD for $\Gamma$-equivariant systems}}
\end{minipage}
\begin{minipage}{0.25\linewidth}
\begin{itemize}
\item{Pick a dictionary of $N$ observables}
\item{Evaluate the observables at data points $x_i$ and $y_i$}
\item{Evaluate the entries of $G,A$: \color{red}{$N^2$ elements}}
\item{Obtain $G^{+}$: \color{red}{$N\times N$ matrix}}
\item{Find $K=G^{+} A$: \color{red}{$N\times N$ matrices}}
\item{Find the eigendecomposition of $K$: \color{red}{$N\times N$ matrix}}
\end{itemize}
%\vfill\null
\end{minipage}
\begin{minipage}{0.73\linewidth}
$~$
\begin{itemize}
\item{Pick a dictionary of $N$ observables}
\item{Identify the symmetry $\Gamma$ of the system, find the irreducible representations of $\Gamma$}
\item{Change the basis to a $\Gamma$-symmetric basis using \cref{eq:projection operator} and \cref{eq:basis} : \color{red}{multiplying at most $N/|\Gamma|$ $|\Gamma|\times|\Gamma|$ matrices by vectors $|\Gamma|\times 1$.}\color{black}{ Let $N_p$ be the number of functions obtained from applying a projection operator $\mathcal{P}^{p}$ corresponding to $p$ th irreducible representation of $\Gamma$ (e.g., $N_p=N/|\Gamma|$ for cyclic groups).}}
\item{Evaluate the observables at data points $x_i$ and $y_i$, add trajectories to reflect the symmetries if necessary}
\item{To obtain the blocks $K_{pq}$ of $K$ (each isotypic component corresponds to $d_p$ blocks), for each $p$}:
\begin{itemize}
\item{Evaluate the entries of $G_{p1},A_{p1}$: \color{red}{$(N_p/d_p)^2$ elements}}
\item{Obtain $G_{p1}^{+}$: \color{red}{$(N_p/d_p)\times (N_p/d_p)$ matrix}}
\item{Find $K_{p1}=G_{p1}^{+} A_{p1}$: \color{red}{$(N_p/d_p)\times (N_p/d_p)$ matrices}}
\item{Find the eigendecomposition of $K_{p1}$: \color{red}{$(N_p/d_p)\times (N_p/d_p)$ matrix}}
\item{The other $K_{pq}$ blocks equal to $K_{p1}$}
\end{itemize}
\item{$K=\bigoplus\limits_{p}\bigoplus\limits_{q=1}^{d_p} K_{pq}$. Its eigenvalues are the eigenvalues of $K_p$, and its eigenvectors only have $N_p$ nonzero elements. Mathematically, eigenvectors $v_{kl}$ of $K$ are of the form $(v_{kl})_i=\bigoplus\limits_p \delta_{pk} v_{pl}$.}
\end{itemize}
%\vfill\null
\end{minipage}
}}%
\caption{EDMD vs modified EDMD \label{ttable} for $\Gamma$-equivariant systems. $|\Gamma|$ is the order of $\Gamma$. The irreducible representations of $\Gamma$ are indexed by $p$ and are $d_p$-dimensional.}
\end{table*}
Koopman eigenfunctions and eigenmodes have many applications in dimensionality reduction, finding the basins of attraction, characterizing coherency between oscillatory systems, etc. Block diagonalizing $K$ allows the efficient computation of the Koopman eigenvalues, eigenfunctions, and modes. 

The kernel DMD is closely related to the EDMD algorithm. It relies on the calculating the eigentriples associated with $K$ from a dual matrix $\hat{K}$ evaluated using a kernel trick commonly applied in machine learning \cite{williams2014kernel}. This method can be computationally advantageous for cases when the number of basis functions exceeds the number of available measurements of the state of the system. We find that the kernel DMD can also be modified to include symmetry considerations in order to optimize the calculations. The method is provided in  \cref{KDMD}.

\subsection{Consequences of symmetry assumptions in the basis}

Assume the data is symmetric as defined by \cref{eq:data symm} with respect to the symmetry group $\Gamma$. A ``perfect" basis is the one respecting the isotypic decomposition of $\Gamma$. Suppose the basis functions belong to isotypic components of $\Sigma\neq\Gamma$. That choice will affect the structure of $K$. We study that structure by evaluating the elements of $A$, since $K$ and $G^{+}$ have the same structure as $A$. 

If the system is $\Gamma$-equivariant and $\Sigma\subset\Gamma$ and the set of actions of $\Sigma$ is a subset of actions of $\Gamma$, the system is also $\Sigma$-equivariant. Thus, picking a basis respecting the isotypic decomposition of $\Sigma$ will have the block diagonal structure corresponding to $\Sigma$. This means that choice of basis results in block diagonal $K$, but its structure does not provide any additional information about the symmetries of the system. 

If the system is $\Gamma$-equivariant and $\Gamma\subset\Sigma$, functions belonging to particular isotypic components of $\Sigma$ are not preserved by the action of $K$. In the case of symmetric trajectories, that can provide information on what the true symmetries of the system are. 

A simple case corresponds to $\Sigma=\Sigma_0\otimes\Gamma$. In this case, every action of $\Sigma_0$ commutes with every action of $\Gamma$. Each isotypic component of $\mathcal{F}$ with respect to $\Sigma$ can be expressed as $\mathcal{F}_{pq}=\mathcal{F}_{\Sigma_0}^p\cap(\mathcal{F}_{\Gamma})^q$, where $\mathcal{F}_{\Sigma_0}^p$ denotes the $p$\ts{th} isotypic component of $\Sigma_0$. In this case, the off-diagonal blocks corresponding to interactions between isotypic components $\mathcal{F}_{p_1 q_1}$ and $\mathcal{F}_{p_2 q_2}$ are zero if $q_1=q_2$, and generally nonzero otherwise. For instance, if a network of three Duffing oscillators similar to the one discussed in, e.g., \cref{ex:iso Duff}, has no permutation symmetry and $\Sigma=Z_2\otimes D_3$, with the action of $Z_2$ being a sign flip in nodal dynamics, the isotypic components corresponding to these $Z_2$ symmetries will not interact, resulting in two blocks in $K$. 

Next, we consider a more general case. We denote the $p$\ts{th} isotypic component of $\mathcal{F}$ with respect to the symmetry group $\Sigma$ by $\mathcal{F}^p_{\Sigma}$. We note that if $\mathcal{F}^p_{\Sigma}\cap\mathcal{F}^{q_1}_{\Gamma}\neq \varnothing$ and $\mathcal{F}^p_{\Sigma}\cap\mathcal{F}^{q_2}_{\Gamma}\neq \varnothing$ where $q_1$ and $q_2$ index different isotypic components of $\Gamma$, then the off-diagonal blocks of $K$ corresponding to interactions between those components are generally nonzero. The condition for $\mathcal{F}^p_{\Sigma}\cap\mathcal{F}^{q}_{\Gamma}\neq \varnothing$ is equivalent to:
\begin{align}
\label{inter}
    \mathcal{P}^p_{\Sigma}\circ(\mathcal{P}^q_{\Gamma}\circ f)\neq 0,
\end{align}
where $f$ is an arbitrary function, and $\mathcal{P}^p_{\Sigma}$ denotes the projection operator onto the $p$\ts{th} isotypic component with respect to the symmetry group $\Sigma$. 
\begin{align}
\label{symm_type}
\begin{gathered}
    \mathcal{P}^p_{\Sigma}\circ(\mathcal{P}^q_{\Gamma}\circ f)\\=
    \sum\limits_{\sigma\in\Sigma}\chi_p(\sigma)^*\sigma_{\rho}\circ \sum\limits_{\gamma\in\Gamma}\chi_q(\gamma)^*\gamma_{\rho}\circ f\\=
    \sum\limits_{\sigma\in\Sigma,\gamma\in\Gamma}\chi_p(\sigma)^*\chi_q(\gamma)^*(\sigma_{\rho}\gamma_{\rho})\circ f .
\end{gathered}
\end{align}
Let $H$ be the set of \textit{left cosets} of $\Gamma$ in $\Sigma$ (defined as $H=\Sigma/\Gamma=\{\sigma\Gamma:\sigma\in\Sigma\}$, where $\sigma\Gamma=\{\sigma\gamma:\gamma\in\Gamma\}$ \cite{golubitsky2003symmetry}). Thus, the condition of \cref{inter} holds if for all $h\in H$:
\begin{align}
\label{char}
\begin{gathered}
    \sum\limits_{\gamma\in\Gamma}\chi^*_{p}(h\gamma^{-1})\chi^*_q(\gamma)=0
\end{gathered}
\end{align}

Using \cref{char}, the structure of $\Gamma$ can be determined given the structure of $K$ and $\Sigma$ used in the calculation. Characters of irreducible representations are available for small order symmetry groups, and scaling up to larger order is possible using computational group theory software. Below is an example for the subgroups of a dihedral group $D_3$.
\begin{example}
Coupled Duffing oscillators: $Z_2\otimes Z_{2,3}$-equivariant system with $Z_2\otimes D_3$ basis functions.

We consider different coupling schemes of networks of 3 Duffing oscillators shown in \cref{subfig:Z2_Z3} and \ref{subfig:Z2_Z2}. We first note that the $Z_2$ symmetry generated by a sign flip is still present in the system for both cases, so two non-interacting blocks corresponding to irreducible representations of that group with respect to that action are still present. Now we focus on the structure of $K$ within each of these non-interacting blocks. 

First, let the function dictionary symmetry be $\Sigma=D_3=\langle r,\kappa | r^3=\kappa^2=e, \kappa r \kappa = r^{-1} \rangle$ and the true symmetry of the system be $\Gamma=Z_3=\langle r| r^3=e\rangle$, where $r_{D_3}$ and $r_{Z_3}$ have the same action. The isotypic component decomposition of $D_3$ is defined in \cref{ex:iso Duff} and can be written as $\mathcal{F}=\mathcal{F}_{tr,D_3}\oplus\mathcal{F}_{sign,D_3}\oplus\mathcal{F}_{st,D_3}$. The isotypic component decomposition of $Z_3$ is defined as $\mathcal{F}=\mathcal{F}_{tr,Z_3}\oplus\mathcal{F}_{\omega,Z_3}\oplus\mathcal{F}_{\omega^2,Z_3}$ ($Z_3$ has 3 1-dimensional irreducible representations with $\chi_{tr}(r)=1$, $\chi_{\omega}(r)=\omega$, $\chi_{\omega^2}(r)=\omega^2$). We note that:
\begin{itemize}
    \item $\mathcal{F}_{tr,Z_3}\cap \mathcal{F}_{tr,D_3}\neq\varnothing$
    \item $\mathcal{F}_{tr,Z_3}\cap \mathcal{F}_{sign,D_3}\neq\varnothing$
    \item $\mathcal{F}_{\omega,Z_3}\cup\mathcal{F}_{\omega^2,Z_3}=\mathcal{F}_{st,D_3}$
\end{itemize}
Thus, the off block-diagonal structure of $K$ is defined by:
\begin{itemize}
    \item $K\mathcal{F}_{tr,D_3}\cap\mathcal {F}_{sign,D_3}\neq\varnothing$
    \item $K\mathcal{F}_{sign,D_3}\cap\mathcal {F}_{tr,D_3}\neq\varnothing$
    \item other off-diagonal blocks are zeros
\end{itemize} 
This structure is illustrated on \cref{subfig: b} and differs from that on  \cref{subfig: a}.

Now, let $\Sigma=D_3$ and $\Gamma=Z_2$. Here, $Z_2=\langle e,\kappa| \kappa^2=e \rangle$, and $\kappa_{D_3}$ and $\kappa_{Z_2}$ have the same action. The isotypic component decomposition of $Z_2$ is defined as $\mathcal{F}=\mathcal{F}_{tr,Z_2}\oplus\mathcal{F}_{sign,Z_2}$.

 We note that:
\begin{itemize}
    \item $\mathcal{F}_{tr,Z_2}\cap \mathcal{F}_{tr,D_3}\neq\varnothing$
    \item $\mathcal{F}_{tr,Z_2}\cap \mathcal{F}_{st,D_3}\neq\varnothing$
    \item $\mathcal{F}_{sign,Z_2}\cap \mathcal{F}_{sign,D_3}\neq\varnothing$
    \item $\mathcal{F}_{sign,Z_2}\cap \mathcal{F}_{st,D_3}\neq\varnothing$
\end{itemize}
Additionally:
\begin{itemize}
    \item $\mathcal{F}_{tr,Z_2}\cap \mathcal{F}_{sign,D_3}=\varnothing$
    \item $\mathcal{F}_{sign,Z_2}\cap \mathcal{F}_{tr,D_3}=\varnothing$
\end{itemize}
Thus, the off block-diagonal structure of $K$ is defined by:
\begin{itemize}
    \item $K\mathcal{F}_{tr,D_3}\cap \mathcal{F}_{sign,D_3}=\varnothing$
    \item other off-diagonal blocks corresponding to interactions between node permutation isotypic components are generally nonzero
\end{itemize} 
This structure is illustrated on  \cref{subfig: c} and differs from that on \cref{subfig: a} and \cref{subfig: b}.
\begin{figure*} 
\begin{subfigure}[b]{.3\textwidth} 
\includegraphics[width=0.9\textwidth]{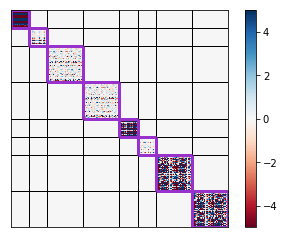}
\caption{$K$ if $\Sigma=\Gamma=D_3$}
\label{subfig: a}
\end{subfigure}
\begin{subfigure}[b]{.3\textwidth}
\includegraphics[width=0.9\textwidth]{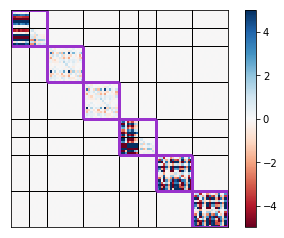}
\caption{$K$ if $\Sigma=D_3$, $\Gamma=Z_3$}
\label{subfig: b}
\end{subfigure}
\begin{subfigure}[b]{.3\textwidth}
\includegraphics[width=0.9\textwidth]{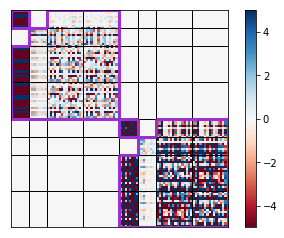}
\caption{$K$ if $\Sigma=D_3$, $\Gamma=Z_2$}
\label{subfig: c}
\end{subfigure}
\caption{Structure of $K$ with basis functions belonging to the isotypic components of $\Sigma=D_3$ for different underlying symmetries of the system}
\label{fig:wrongK}
\end{figure*} 

This example shows that the structure of the approximation $K$ with maximal symmetries assumed provided information about the actual underlying symmetries of the system.

\end{example}

\subsection{Towards realistic systems}
In this manuscript, we provide a general approach for dimensionality reduction in the calculation of Koopman operator approximations by exploiting the underlying symmetries of the system’s dynamics and structure. The exact scaling achieved by the reduction depends on the structure of the symmetry group of the dynamical system, specifically, the number of irreducible representations of the symmetry group and their dimensionality. 

The results outlined in this manuscript, similar to most of the other literature related to dynamical systems with symmetries, are immediately applicable in the case of existence of exact symmetries in nonlinear dynamics. That is the case when the system is completely deterministic and the initial conditions respect the symmetries of the system. If the symmetries of the system are known and the available trajectories are deterministic, it is always possible to reconstruct the trajectories that are related via the symmetry group of the system. 
Then, a full set of trajectories respecting the symmetries of the system can be used to approximate the Koopman operator and its eigendecomposition.

However, in many systems that information is not necessarily available ahead of time and the symmetries are not present in data, even if the initial conditions are symmetric, because of the presence of noise in the system.  Some of the examples of not fully symmetric data include the following cases and their combinations:

\begin{itemize}
\item Deterministic systems with measurement noise.
DMD for systems with measurement noise and possible ways to correct for it are presented in Ref. \cite{dawson2016characterizing}. It is shown in \cref{app:sensor} that in this case the expected values of off-diagonal elements of $K$ computed using the EDMD are zero, so the block decomposition may still be applicable.
\item{Stochastic systems with symmetric initial conditions and process noise.}
DMD applied to the systems with process noise is studied, for instance, in Ref.
\cite{bagheri2014effects}.
\item{Systems with imperfect symmetries due to sampling and unknown underlying symmetries.}
\item{Systems with imperfect symmetries in dynamics (e.g. slight parameter mismatch).}
\end{itemize}

All these cases require separate treatment, and whether the isotypic component decomposition is still useful in computing the Koopman operator approximation will vary depending on specific characteristics of the data available from the system, such as the strength of noise or the trajectory sampling characteristics.  

\section{Conclusion}
In this manuscript, we apply tools from group theory and representation theory to study the structure of the Koopman operator for equivariant dynamical systems. This approach can be applied to systems with permutation symmetries (e.g. networks symmetric under node permutations where the information about the symmetries is contained in the adjacency matrix), systems with intrinsic dynamical symmetries, and systems with both types of symmetries present. We find that the operator itself and its approximations can be block diagonalized using a symmetry basis that respects the isotypic component structure related to the underlying symmetry group and the actions of its elements. For the approximation matrix to be exactly block diagonal, the data must respect the symmetries of the system. That can be readily accomplished if the underlying symmetry is known ahead of time (e.g., the topology of the network is known). Symmetry considerations are applicable to both EDMD and kernel DMD, which means they become useful both in the regime when the number of observables is much larger than that of measurements and vice versa.

% If in two-column mode, this environment will change to single-column format so that long equations can be displayed. 
%\begin{widetext}
%$$\mbox{put long equation here}$$
%\end{widetext}

% Figures should be put into the text as floats. 
% Use the graphics or graphicx packages (distributed with LaTeX2e).
% Put the label that you will use with \ref{} command in the braces of the \label{} command.
%
% \begin{figure}
% \end{figure}

% \begin{table}
% \end{tabular}
% \end{table}

%\begin{acknowledgments}
% Put your acknowledgments here.
%\end{acknowledgments}

\begin{acknowledgments}
This work is supported by the U.S. Army
Research Laboratory and the U.S. Army Research Office under contract number
W911NF-13-1-0340.  The authors thank Jordan Snyder, Mehran Mesbahi, and Afshin Mesbahi for useful discussions.
\end{acknowledgments}

\appendix
\section{Block diagonalization of isotypic components obtained from $d_p$-dimensional irreducible representations} \label{app:iso}

We show that $d$-dimensional irreducible representations of $\Gamma$ yield identical blocks of $\mathcal{K}$ in the isotypic component basis obtained using the unitary irreducible representations of $\Gamma$.

Let the function space be decomposed into isotypic components according to the actions of the symmetry group $\Gamma$ of order $|\Gamma|$, $\gamma_{\rho}\in\Gamma_{\rho}$: $\mathcal{F}=\mathcal{F}_1\oplus...\oplus \mathcal{F}_N$, where $N$ is the number of irreducible representations of $\Gamma$. Let $\mathcal{F}_p$ be one of these isotypic components with a corresponding \textit{unitary} irreducible representation with elements $R_p(\gamma)$ corresponding to $\gamma\in\Gamma$, and let $d_p$ be the dimensionality of that representation. 

The projection operator is defined as:
\begin{align}
\mathcal{P}^p_{mn}=\dfrac{d_p}{|\Gamma|}\sum\limits_{\gamma\in\Gamma} [R_p(\gamma)]^*_{mn} \gamma_{\rho}.
\end{align}
It acts on $f\in\mathcal{F}$ to produce sets of projected functions according to:
\begin{align}
\xi_{mn}^p=\mathcal{P}^p_{mn} \circ f.
\end{align}
We already know that $\mathcal{K} \xi^p_{mn} = h_p$, where $h_p\in\mathcal{F}_p$. The subspace $\mathcal{F}_p$ can be decomposed into $d_p$ components $\mathcal{F}_p=\mathcal{F}_{p,1}\oplus...\oplus \mathcal{F}_{p,d_p}$, where $\mathcal{F}_{p,m} =\{g|g=\mathcal{P}_{mn}\circ f,~f\in\mathcal{F}, n=1,...,d_p\}$. This is a well-defined decomposition since $\langle \mathcal{P}^p_{mn} f, \mathcal{P}^p_{kl}h \rangle = \langle f, \mathcal{P}^p_{nm}\mathcal{P}^p_{kl}h \rangle = \langle f, \delta_{mk}\mathcal{P}^p_{nl}h \rangle$ \cite{cornwell1997group} can be nonzero only when $m=k$. 

We want to show that $\mathcal{K} \xi^p_{mn}\in \mathcal{F}_{p,m}$ (also true for any linear operator that commutes with the action of the symmetry group). Since the operator commutes with the actions of the group:
\begin{align}
\mathcal{K}\mathcal{P}^p_{mn} \circ f = \mathcal{P}^p_{mn} \mathcal{K} \circ f = \mathcal{P}^p_{mn} \circ h \in \mathcal{F}_{p,m}.
\end{align}
Here, $\mathcal{K} \circ f \equiv h$.

Let $f_{\Gamma}=\{\gamma\circ f|\gamma\in\Gamma\}$. Any set of linearly independent functions that span $f_{\Gamma}$ can be transformed into a symmetry respecting basis obtained by calculating all the projections $\mathcal{P}^p_{mn}\circ f_{\gamma}$, where $f_{\gamma}\in f_{\Gamma}$. That corresponds to a block diagonal form of the Koopman operator.

We've already shown that $K$, the approximation of $\mathcal{K}$, also commutes with the actions of the elements of $\Gamma$ for $\Gamma$-equivariant dynamical systems with $\Gamma$-equivariant data. Thus, we can obtain an observable dictionary that block diagonalizes $K$ into $|\Gamma|$ blocks, where each $d_p$-dimensional irreducible representation results in $d_p$ $d_p\times d_p$-dimensional blocks.

Additionally, suppose $\mathcal{K} \mathcal{P}^p_{mn} \circ f = h$, then $\mathcal{K} \mathcal{P}^p_{kn} \circ f = \mathcal{P}^p_{km}\mathcal{K}\mathcal{P}^p_{mn} \circ f= \mathcal{P}^p_{km}\circ h$. This gives us the relation between blocks in $\mathcal{K}$ corresponding to the same irreducible representation $p$. In context of the approximation $K$, it means that we get that $K_{p,i}$ (blocks corresponding to $\psi\in\mathcal{F}_{p,i}$) are equal for all $i$ (for data respecting the symmetries of the system and a proper ordering of basis functions).

\section{Commutativity of $K$ and $\gamma_{\rho}$ acting in function space} \label{app: K commutes}

We show that $G$ and $G^{+}$ ($^{+}$ denotes the Moore-Penrose pseudoinverse) commute. We note that $G$ is a Hermitian matrix since:
\begin{align}
\begin{gathered}
    G^{*}=\left(\sum\limits_m \Psi^*(x_m)\Psi(x_m)\right)^*\\
    =\sum\limits_m\Psi^*(x_m)\Psi(x_m)=G
\end{gathered}
\end{align}
Thus, $G$ is also normal, i.e. $GG^{*}=G^{*}G$. We show that if $G$ is normal, $GG^{+}=G^{+}G$.

Two of the criteria that define the Moore-Penrose pseudoinverse \cite{penrose1955generalized} state that $G^{+}=G^+ G G^+$ and $(GG^+)^* = GG^+$. It follows that the following relation holds: $G^+=G^+(GG^+)^*=G^+(G^+)^* G^*$. Using that relation \cite{ben2003generalized} and commutativity of $^+$ and $^*$ operations, we obtain:
\begin{align}
\begin{gathered}
G^+ G = G^+ (G^+)^* G^* G = (G^+)^* G^+ G G^*\\ = (G^+)^* G^*
= (G^+ (G^+)^* G^*)^* G^*\\ = GG^{+} (G^{+})^* G^* = GG^+.
\end{gathered}
\end{align}

Since the action of $\gamma$ commutes with $A$ and $G$, and since $G$ commutes with $G^{+}$, the action of $\gamma$ commutes with $K=G^+A$. which is a Koopman operator approximation.

\section{Change of basis and the EDMD approximation} \label{app: basis}
We show that rotating the observable dictionary preserves the symmetries of the reconstructed trajectories.

Suppose we have a basis consisting of dictionary functions $\mathcal{D}_{\psi}$ and a dictionary $\mathcal{D}_{\xi}$ obtained by $\Xi=T\Psi$. Let $\Psi(t)=(\psi_1(x(t))~~\psi_N(x(t)))^T$ and $\Xi(t)=(\xi_1(x(t))~~\xi_N(x(t)))^T$. We show that rotating the dictionary function vector does not affect the trajectory reconstruction: 
\begin{align}
\begin{gathered}
    \Psi_{t+1}= K_{\psi}\Psi_{t}\\
    \Xi_t = T\Psi_{t} \\
    \Xi_{t+1}=  K_{\psi}T\Psi_{t} =  T K_{\psi}\Psi_{t}   = T\Psi_{t+1}
\end{gathered}
\end{align}

Next, we show that the state reconstruction preserves the symmetries of the system. Let $P$ be the action of the symmetry group on the basis functions $\Psi$. We aim to show that if $\Psi_{t+1}=K\Psi_{t}$, then $ P\Psi_{t+1}=K P\Psi_{t}$. It follows directly from the fact that $K$ and $P$ commute:
\begin{align}
    P \Psi_{t+1} = PK \Psi_{t} = K P\Psi_{t}
\end{align}
Thus, the trajectories of basis functions reconstructed using the EDMD approximation are $\Gamma$-equivariant, just like the original system. In particular, this is true in case of the evolution of the full state observable.

\section{Applicability to kernel methods}{\label{KDMD}}

Kernel DMD introduced in Ref. \cite{williams2014kernel} is a variant of approximating the Koopman operator matrix most efficient when the number of measurement points is much smaller than the number of basis functions. Kernel DMD relies on evaluating $\hat{G}$ and $\hat{A}$ using the kernel method. Their elements can be found by indirectly evaluating the inner products in the basis function space: $k(x_m,y_n)=\Psi(x_m)\Psi(y_n)^*$ (e.g., if $k$ is a polynomial kernel, $k(x,y)=(1+xy^T)^{\alpha}$). We note that $k(\gamma x, \gamma y)= k(x,y)$ due to the properties of inner products.

In kernel DMD, $\hat{G}_{ij}=k(x_i,x_j)$ and $\hat{A}_{ij}=k(x_i,y_j)$. The eigendecomposition of $\hat{G}=Q\Sigma^2 Q^T$ is then used to find the matrix $\hat{K}$ and use it in computing the eigendecomposition of the Koopman operator approximation matrix $K$:
\begin{align}
    \hat{K}=(\Sigma^{+}Q^T)\hat{A}(Q\Sigma^{+}).
\end{align}

Again, we pick a particular order of group elements similarly to \cref{eq: order}:
\begin{align} \label{eq: orderr}
\begin{gathered}
    \Psi_x=\begin{pmatrix}
        \Psi(\gamma_1 \textbf{x})\\
        \vdots\\
        \Psi(\gamma_{|\Gamma|} \textbf{x})\\
    \end{pmatrix},\\
    \text{where}~~\Psi(\textbf{x})=\begin{pmatrix}
    \Psi_1(x_1)&...&
    \Psi_N(x_1)\\
    \vdots&\ddots&\vdots\\
    \Psi_1(x_{M/|\Gamma|})&...&
    \Psi_N(x_{M/|\Gamma|})
    \end{pmatrix}
\end{gathered}
\end{align}
We also construct a permutation representation of the group with elements denoted by $P_{\gamma_i}$ as defined in \cref{eq:permutation}.

By Cayley's theorem, such permutations form a group isomorphic to $\Gamma$. Determining the actions $P_{\gamma_i}$ of the group generators is sufficient to find the actions of all the group elements. Let $\bm{P}_{\gamma_k}=P_{\gamma_k}\otimes I_{n\times n}$. We note that $\left( \bm{P}_{\gamma_k} \right)^*=(\bm{P}_{\gamma_k})^{-1}$. It can be shown that:
\begin{align}\label{eq:comm}
    \bm{P}_{\gamma_i} \hat{G} = \hat{G} \bm{P}_{\gamma_i},~~\bm{P}_{\gamma_i} \hat{A} = \hat{A} \bm{P}_{\gamma_i}.
\end{align}
We do so for $\hat{A}$, and the proof for $\hat{G}$ is equivalent. We find that:
\begin{align}
\begin{gathered}
(\bm{P}_{\gamma_i} \hat{A})_{kl}=\hat{A}_{pl}=k(x_p,y_l),~~\gamma_p=\gamma_i\gamma_k,\\
(\hat{A} \bm{P}_{\gamma_i})_{kl}=\hat{A}_{kq}=k(x_k,y_q),~~\gamma_q=\gamma_i^{-1}\gamma_l.
\end{gathered}
\end{align}
And finally, $k(x_p,y_l)=k(\gamma_i x_k, \gamma_i y_q)=k(x_k, y_q)$.

Since the relation \ref{eq:comm} holds, the same reasoning can be applied to block diagonalize the matrix $\hat{K}$. It is sufficient to apply the projection operator \cite{stiefel2012group}:
\begin{align} \label{eq:projection operator2}
\mathcal{P}^p_{mn}=\dfrac{d_p}{|\Gamma|}\sum\limits_{\gamma\in\Gamma} [R_p(\gamma)]^*_{mn} P_{\gamma_i}.
\end{align}
This projection operator is analogous to the one introduced in equation \ref{eq:projection operator1}, except the symmetry group in this case acts by permuting the group elements.

We can apply the singular value decomposition of $\mathcal{P}$ to obtain the basis for the projection subspaces of irreducible representations (isotypic components).
We form the transformation matrix $T$ by finding the singular value decomposition (SVD) and stacking its eigenvectors as rows of $T$ such that $\bm{T}=T\otimes I_{n\times n}$.

\begin{table*} [t]
\fbox{%
    \parbox{1.\textwidth}{
\begin{minipage}{0.25\linewidth} 
\centering
{\textbf{Standard kernel DMD}} 
\end{minipage}
\begin{minipage}{0.73\linewidth}
\centering
{\textbf{Kernel DMD for $\Gamma$-equivariant systems}}
\end{minipage}
\begin{minipage}{0.25\linewidth}
\begin{itemize}
\item{Pick a dictionary of $N$ observables}
\item{Evaluate the kernel functions at data points $x_i$ and $y_i$}
\item{Evaluate the entries of $\hat{G},\hat{A}$: \color{red}{$M^2$ elements}}
\item{Obtain $\hat{G}^{+}$: \color{red}{$M\times M$ matrix}}
\item{Find\\ $    \hat{K}=(\Sigma^{+}Q^T)\hat{A}(Q\Sigma^{+})
$: \color{red}{$M\times M$ matrices}}
\item{Find the eigendecomposition of $\hat{K}$: \color{red}{$M\times M$ matrix}}
\end{itemize}
%\vfill\null
\end{minipage}
\begin{minipage}{0.73\linewidth}
$~$
\begin{itemize}
\item{Pick a dictionary of $N$ observables}
\item{Identify the symmetry $\Gamma$ of the system, find the irreducible representations of $\Gamma$}
\item{Change the basis to a $\Gamma$-symmetric basis using \cref{eq:permutation} and \cref{eq:projection operator2}}
\item{Evaluate the observables at data points $x_i$ and $y_i$, add trajectories to reflect the symmetries if necessary}
\item{To obtain the blocks $\hat{K}_{pq}$ of $\hat{K}$ (each isotypic component corresponds to $d_p$ blocks), for each $p$}:
\begin{itemize}
\item{Evaluate the entries of $\hat{G}_{p1},\hat{A}_{p1}$: \color{red}{$(M_{p}/d_p)^2$ elements}}
\item{Obtain $\hat{G}_{p1}^{+}$: \color{red}{$(M_{p}/d_p)\times (M_{p}/d_p)$ matrix}}
\item{Find $\hat{K}_{p1}=\hat{G}_{p1}^{+} \hat{A}_{p1}$: \color{red}{$(M_{p}/d_p)\times (M_{p}/d_p)$ matrices}}
\item{Find the eigendecomposition of $\hat{K}_{p_1}$: \color{red}{$(M_{p}/d_p)\times (M_{p}/d_p)$ matrix}}
\item{The other $\hat{K}_{pq}$ blocks equal to $\hat{K}_{p1}$}
\end{itemize}
\item{$\hat{K}=\bigoplus\limits_{p}\bigoplus\limits_{q=1}^{d_p} \hat{K}_{pq}$. Its eigenvalues are the eigenvalues of $\hat{K}_p$, and its eigenvectors only have $M_p$ nonzero elements. Mathematically, eigenvectors $v_{kl}$ of $K$ are of the form $(v_{kl})_i=\bigoplus\limits_p \delta_{pk} v_{pl}$.}
\end{itemize}
%\vfill\null
\end{minipage}
}}
\caption{kernel DMD vs modified kernel DMD for $\Gamma$-equivariant systems. $|\Gamma|$ is the order of $\Gamma$. The irreducible representations of $\Gamma$ are indexed by $p$ and are $d_p$-dimensional. Here, $M$ be the number of data points used by the algorithm, and $\{(x_m,y_m)\}$ respect the symmetries of the system.}
\label{tttable}
\end{table*}

Similarly to EDMD, the isotypic component basis simplifies calculating the approximations of $\hat{K}$.
\begin{align}
\begin{gathered}
   \hat{A}_D=\bigoplus\limits_p \bigoplus\limits_q \hat{A}_{pq},\\
   \hat{G}_D=\bigoplus\limits_p \bigoplus\limits_q \hat{G}_{pq}=\bigoplus\limits_p \bigoplus\limits_q Q_{pq}\Sigma^2_{pq}Q_{pq}^T,\\
   \hat{K}_D=
    \bigoplus\limits_p \bigoplus\limits_q(\Sigma_{pq}^{+}Q_{pq}^T)\hat{A}_{pq}(Q_{pq}\Sigma_{pq}^{+}).
\end{gathered}
\end{align}
The modification is summarized in table \ref{tttable}. 

Finally, the approximations of Koopman eigenvalues, eigenfunctions, and eigenmodes can be calculated using $K_D$, as shown in Ref. \cite{williams2014kernel}.

\section{Deterministic systems with sensor noise} \label{app:sensor}

Transfer operators with process and measurement noise were also studied in Ref. \cite{sinha2018robust}. Characterizing and correcting for the effect of sensor noise in DMD is discussed in Ref. \cite{dawson2016characterizing}. We need to extend the results to EDMD to quantify the effect of sensor noise on the structure of the matrix $K$. The main modification that needs to be made is the consideration of the effect of the noise in measuring $X$ and $Y$ on $\Psi_x$ and $\Psi_y$.

Let $X$ and $Y$ be matrices analogous to $\Psi_x$ and $\Psi_y$ corresponding to the full-state observable evaluated at discrete time steps. We denote the sensor noise matrices by $N_x$ and $N_y$, so that the measured $X_n$ and $Y_n$ can be found from $X_n=X+N_x$ and $Y_n=Y+N_y$. We assume that the noise distributions respect the symmetries of the system, which might be the case, for instance, for symmetric networks. Moreover, we assume that the noise is state-independent.

We can form vectors $\Psi_{xn}$ and $\Psi_{yn}$ that can be used to find the approximation $K$ using EDMD:
\begin{align}
K_n=\Psi_{xn}^{+}\Psi_{yn}.
\end{align}
Here, $(\Psi_{xn})_{ij}=\psi_j((X_n)_i)$, $(\Psi_{yn})_{ij}=\psi_j((Y_n)_i)$, and $N_{\Psi,x}$ and $N_{\Psi,y}$ correspond to noise matrices obtained as:
\begin{align}
\begin{gathered}
(N_{\Psi,x})_{ij}=\psi_j(X_i+N_{x,i})-\psi_j(X_i).
\end{gathered}
\end{align}

We aim to show that $\mathbb{E}(PK_n)=\mathbb{E}(K_nP)$, meaning that the expected value of the Koopman operator $K_n$ commutes with the permutation matrix corresponding to an element of the symmetry group. If that is the case, then the expected values of the off-block-diagonal elements of $K_n$ in a symmetry adapted basis as defined in \cref{eq:projection operator} are zero. To do that, we can express $K_n$ as:
\begin{align} \label{Eq: noisy K}
\begin{gathered}
K_n=\Psi_{xn}^{+}\Psi_{yn}=(\Psi_{x}+N_{\Psi,x})^{+}(\Psi_{y}+N_{\Psi,y})\\=
((\Psi_{x}+N_{\Psi,x})^*(\Psi_{x}+N_{\Psi,x}))^{+}
\\ (\Psi_{x}+N_{\Psi,x})^{*}(\Psi_{y}+N_{\Psi,y})\\=(\Psi_{x}^{*}\Psi_{x}+\Psi_{x}^{*}N_{\Psi,x}+N_{\Psi,x}^{*}\Psi_{x}+N_{\Psi,x}^{*}N_{\Psi,x})^{+}\\(\Psi_{x}^{*}\Psi_{y}+\Psi_{x}^{*}N_{\Psi,y}+N_{\Psi,x}^{*}\Psi_{y}+N_{\Psi,x}^{*}N_{\Psi,y}).
\end{gathered}
\end{align}
If the inverse of the first term exists, it can be expanded into the Taylor series with terms of the form below in a weak noise limit. We need to show that:
\begin{align}
    \bm{P}_{\gamma_k}\mathbb{E}(M_1^* M_2...M_{n-1}^*M_n)\\=
    \mathbb{E}(M_1^* M_2...M_{n-1}^*M_n)\bm{P}_{\gamma_k}
\end{align}
Here, the matrices $M_i$ are selected from $N_{\Psi,x/y}$ and $\Psi_{x/y}$.
That follows directly from:
\begin{align}
\begin{gathered}
    \mathbb{E}((M_1\bm{P}_{\gamma_k})^* (M_2\bm{P}_{\gamma_k})...(M_{n-1}\bm{P}_{\gamma_k})^*(M_n\bm{P}_{\gamma_k}))\\
    =\bm{P}^{-1}_{\gamma_k}\mathbb{E}(M_1^* M_2...M_{n-1}^*M_n)\bm{P}_{\gamma_k}\\
    =\mathbb{E}(M_1 M_2^*...M_{n-1}M_n^*)
\end{gathered}
\end{align}
Thus, the expected values of the off-block-diagonal elements of $K_n$ are zero in the isotypic component basis.

%\nocite{*}
\bibliography{aipsamp}
\bibliographystyle{plain}

\end{document}